\documentclass[a4paper,11pt,oneside]{article}

\usepackage{amsmath,amsthm,amsfonts,amssymb,amscd}
\usepackage{graphicx}
\usepackage{mathtools}
\usepackage{hyperref}
\usepackage{pgf,tikz}
\usepackage{mathrsfs}
\usepackage{float}
\usepackage{csvsimple}
\usepackage{enumerate}
\usepackage{booktabs}
\usepackage{pgfplots}
\pgfplotsset{compat=1.15}
\usepackage[noend]{algpseudocode}
\usepackage{graphicx,bm,xcolor}
\usepackage{algorithm}
\usepackage{pdfpages}
\usepackage{algpseudocode}
\usepackage{stmaryrd}
\usepackage{hyperref}
\usetikzlibrary{arrows}
\usepackage{caption}
\usepackage[T1]{fontenc}
\usepackage{t1enc}
\usepackage[polish,german,hungarian,english]{babel}
\usepackage{verbatim} 
\usepackage{xcolor}

\usepackage{url}

\usepackage{siunitx}

\newcommand\plot[1]{\let\frame\relax
	\frame{\includegraphics[clip,trim=0 220 0 60,width=8cm]{#1}}}
\newcommand{\dx}{\text{ d}x}

\usepackage{authblk}

\usepackage[top=2cm, bottom=2cm, left=2cm, right=2cm]{geometry}

\theoremstyle{plain} 
\newtheorem{theorem}{Theorem}[section]

\newtheorem{Remark}[theorem]{Remark}

\theoremstyle{definition} %

\theoremstyle{remark} %

\usepackage{newunicodechar}
\newunicodechar{ℝ}{\mathbb{R}}
\newunicodechar{Ψ}{\Psi}
\newunicodechar{α}{\alpha}
\newunicodechar{∇}{\nabla}
\newunicodechar{σ}{\sigma}
\newunicodechar{ν}{\nu}
\newunicodechar{ρ}{\rho}
\newunicodechar{ϑ}{\vartheta}
\newunicodechar{Γ}{\Gamma}

\usepackage{xcolor}

\usepackage{pgfplots}
\pgfplotsset{compat=1.15}
\usepackage{mathrsfs}
\usetikzlibrary{arrows}
\usepackage{tikz-3dplot}

\newif\ifMAKEPICS
\MAKEPICSfalse 

\ifMAKEPICS
\usepackage[cleanup,subfolder]{gnuplottex}
\usepackage{xparse}

\ExplSyntaxOn
\DeclareExpandableDocumentCommand{\convertlen}{ O{cm} m }
{
	\dim_to_decimal_in_unit:nn { #2 } { 1 #1 } cm
}
\ExplSyntaxOff
\fi

\renewcommand{\dx}{\textrm{ d}x}
\usepackage{graphicx}
\begin{document}

\title{Locally different models in a checkerboard pattern with mesh adaptation and error control for multiple quantities of interest}

\author[1,2]{B. Endtmayer}

\affil[1]{Leibniz Universit\"at Hannover, Institut f\"ur Angewandte
	Mathematik, Welfengarten 1, 30167 Hannover, Germany}

\affil[2]{Cluster of Excellence PhoenixD (Photonics, Optics, and
	Engineering -- Innovation Across Disciplines), Leibniz Universit\"at Hannover, Germany}

\maketitle                   
\begin{abstract}
	In this work, we apply multi-goal oriented error estimation to the finite element method. In particular, we use the dual weighted residual method and apply it to a model problem. This model problem consist of  locally different coercive partial differential equations in a checkerboard pattern, where the solution is continuous across the interface. In addition to the error estimation, the error can be localized using a partition of unity technique. The resulting adaptive algorithm is substantiated with a numerical example.
\end{abstract}
\section{Introduction}
\label{sec: intro}
In many applications partial differential equations (PDEs) have to be solved. One possible method to approximate the solution of these PDEs is the finite element method \cite{courant1943variational}.
A posteriori error estimation and mesh adaptivity are well developed methodologies for the finite element method; see for example \cite{Cia87,braess2013finite,dorfler1996convergent,AinsworthOden:2000,verfurth2013posteriori,babuvska1978posteriori,repin2008posteriori,BeRa01} and the references therein.
However, the solution of the PDEs might not of primary interest, but some quantity of interest, evaluated at the solution. 
Goal oriented error estimation  \cite{BeRa01,BeRa96,KerPruChaLaf2017,EnLaNeWiWo20,BREVIS2021186,BruZhuZwie16,MeiRaVih109,BaRa03,bause2021flexible,bruchhauser2020dual,dolejvsi2021goal,AhEndtSteiWi22,BeuEndtLaWi23,endtmayer2024goal} allows us the estimate the error in the given quantity of interest. This approach can be extended to multiple quantities of interested; see for instance \cite{HaHou03,Ha08,PARDO20101953,BruZhuZwie16,EnLaWiWoNeiPAMM2019,KerPruChaLaf2017,EndtWi17,kergrene2018goal,EnLaWi18,EndtLaRiSchafWi24_book_chapter} and the references therein.
This work is devoted to estimate the error for a combined quantity of interest, which bounds the error of all quantities of interest, for a problem with locally different models, which are arranged in a checkerboard pattern..  In Section~\ref{sec: model problem}, we introduce the model problem with locally different models. Section~\ref{sec: disc} describes the used finite element discretization. Goal oriented error estimation is discussed in Section~\ref{sec: Goafem} and the extension to multiple quantities of interest is given in Section~\ref{sec: multigoal}. Finally the developments are substantiated with a numerical example in Section~\ref{sec: numerical results}.

\section{The Model Problem}
\label{sec: model problem}
\begin{figure}[H]
\begin{minipage}{0.49\linewidth}
	\definecolor{zzttqq}{rgb}{0.6,0.2,0}
	\definecolor{xdxdff}{rgb}{0.5,0.2,1}
	\definecolor{ududff}{rgb}{0.2,0.5,1}
	\definecolor{uuuuuu}{rgb}{0.2,0.2,0.2}
	\begin{tikzpicture}[line cap=round,line join=round,>=triangle 45,x=3.4cm,y=3.4cm]
	\clip(-0.2,-0.2) rectangle (2.2,2.2);
	\fill[line width=2pt,color=uuuuuu,fill=zzttqq,fill opacity=0.00000000149011612] (0,0) -- (2,0) -- (2,2) -- (0,2) -- cycle;
	\fill[line width=1pt,color=zzttqq,fill=zzttqq,fill opacity=0.4] (0,1) -- (1,1) -- (1,2) -- (0,2) -- cycle;
	\fill[line width=1pt,color=zzttqq,fill=xdxdff,fill opacity=0.4] (1,1) -- (2,1) -- (2,2) -- (1,2) -- cycle;
	\fill[line width=1pt,color=zzttqq,fill=ududff,fill opacity=0.4] (0,0) -- (1,0) -- (1,1) -- (0,1) -- cycle;
	\fill[line width=1pt,color=zzttqq,fill=uuuuuu,fill opacity=0.4] (1,0) -- (2,0) -- (2,1) -- (1,1) -- cycle;
	\draw [line width=3pt] (0,2)-- (2,2);
	\draw [line width=3pt] (2,2)-- (2,0);
	\draw [line width=3pt] (2,0)-- (0,0);
	\draw [line width=3pt] (0,0)-- (0,2);
	\draw [line width=2pt,color=green] (0,1)-- (2,1);
	\draw [line width=2pt,color=green] (1,2)-- (1,0);
	\begin{scriptsize}
	\draw[color=uuuuuu] (1.5,1.5) node {\normalsize$-\Delta_\varepsilon^p u$};
	\draw[color=uuuuuu] (0.5,1.5) node {\normalsize$-\Delta u+u^3$};
	\draw[color=uuuuuu] (1.5,0.5) node {\normalsize$-\textrm{div}\left((\frac{1}{10}+\frac{1}{1+u^2})\nabla u \right)$};
	\draw[color=uuuuuu] (0.5,0.5) node {\normalsize$-\Delta u$};
	\draw[color=uuuuuu] (1,2.1) node {\normalsize$\partial \Omega$};
	\draw[color=uuuuuu] (2.1,1) node {\normalsize$\partial \Omega$};
	\draw[color=uuuuuu] (-0.1,1) node {\normalsize$\partial \Omega$};
	\draw[color=uuuuuu] (1,-0.1) node {\normalsize$\partial \Omega$};
	\end{scriptsize}
	\end{tikzpicture}
\end{minipage}%
\begin{minipage}{0.49\linewidth}
	\definecolor{zzttqq}{rgb}{0.6,0.2,0}
	\definecolor{xdxdff}{rgb}{0.5,0.2,1}
	\definecolor{ududff}{rgb}{0.2,0.5,1}
	\definecolor{uuuuuu}{rgb}{0.2,0.2,0.2}
	\begin{tikzpicture}[line cap=round,line join=round,>=triangle 45,x=3.4cm,y=3.4cm]
	\clip(-0.2,-0.2) rectangle (2.2,2.2);
	\fill[line width=2pt,color=uuuuuu,fill=zzttqq,fill opacity=0.00000000149011612] (0,0) -- (2,0) -- (2,2) -- (0,2) -- cycle;
	\fill[line width=1pt,color=zzttqq,fill=zzttqq,fill opacity=0.4] (0,1) -- (1,1) -- (1,2) -- (0,2) -- cycle;
	\fill[line width=1pt,color=zzttqq,fill=xdxdff,fill opacity=0.4] (1,1) -- (2,1) -- (2,2) -- (1,2) -- cycle;
	\fill[line width=1pt,color=zzttqq,fill=ududff,fill opacity=0.4] (0,0) -- (1,0) -- (1,1) -- (0,1) -- cycle;
	\fill[line width=1pt,color=zzttqq,fill=uuuuuu,fill opacity=0.4] (1,0) -- (2,0) -- (2,1) -- (1,1) -- cycle;
	\draw [line width=3pt] (0,2)-- (2,2);
	\draw [line width=3pt] (2,2)-- (2,0);
	\draw [line width=3pt] (2,0)-- (0,0);
	\draw [line width=3pt] (0,0)-- (0,2);
	\draw [line width=2pt,color=green] (0,1)-- (2,1);
	\draw [line width=2pt,color=green] (1,2)-- (1,0);
	\begin{scriptsize}
	\draw[color=uuuuuu] (1.5,1.5) node {\normalsize$\Omega_2$};
	\draw[color=uuuuuu] (0.5,1.5) node {\normalsize$\Omega_1$};
	\draw[color=uuuuuu] (1.5,0.5) node {\normalsize$\Omega_3$};
	\draw[color=uuuuuu] (0.5,0.5) node {\normalsize$\Omega_4$};
	\draw[color=green] (1.1,0.85) node {\Large$I$};
	\draw[color=uuuuuu] (1,2.1) node {\normalsize$\partial \Omega$};
	\draw[color=uuuuuu] (2.1,1) node {\normalsize$\partial \Omega$};
	\draw[color=uuuuuu] (-0.1,1) node {\normalsize$\partial \Omega$};
	\draw[color=uuuuuu] (1,-0.1) node {\normalsize$\partial \Omega$};
	\end{scriptsize}
	\end{tikzpicture}
\end{minipage}%
\caption{The locally different models of the partial differential equations on the domain $\Omega$ (left) and the sub domains $\Omega_i$ for $i \in \{1,2,3,4\}$ with the interface $I$ (right) .\label{fig: config}}
\end{figure}
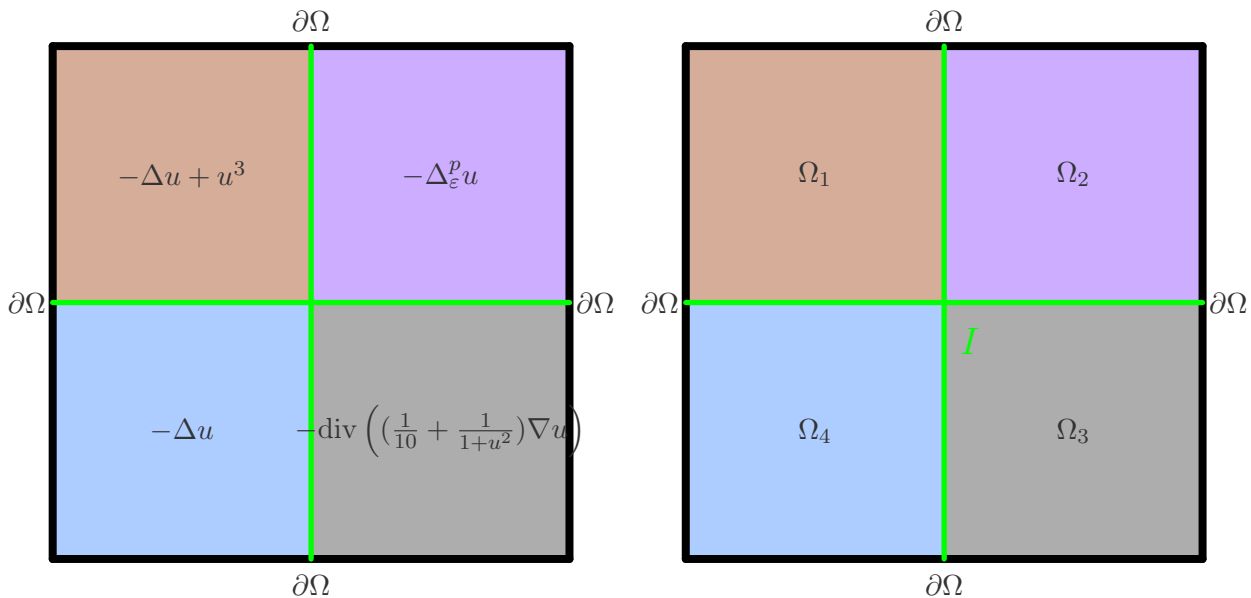

Let $\Omega$ be $(-1,1)^2 \subset \mathbb{R}^2$, $\Omega_1:=\{(x_1,x_2) \in \Omega: x_1<0,x_2>0\}$, $\Omega_2:=\{(x_1,x_2) \in \Omega: x_1>0,x_2>0\}$, , $\Omega_3:=\{(x_1,x_2) \in \Omega: x_1>0,x_2<0\}$ and  $\Omega_4:=\{(x_1,x_2) \in \Omega: x_1<0,x_2<0\}$ as depicted in Figure~\ref{fig: config}. Furthermore, let $W_1^{p(x)}(\Omega)$ be defined as in \cite{mashiyev2010existence,diening2004riesz,kovavcik1991spaces} and let $$p(x):=\begin{cases}
2 \qquad \qquad\qquad \qquad \qquad \textrm{ }\textrm{ in } \Omega_1 \cup \Omega_3 \cup \Omega_4,\\
(1+2x_1)(1+2x_2) \qquad \quad \textrm{in } \Omega_2
\end{cases}.$$
Furthermore, $W_{1,0}^{p(x)}(\Omega):=\{u \in W_1^{p(x)}(\Omega): u=0 \text{ on } \partial \Omega\}$.
The model problem is formally given by:
Find $u \in V:= W_{1,0}^{p(x)}(\Omega)$ such that
\begin{align*}
	{A}(u)&=f\qquad \textrm{in }\Omega, \\
	u&=0 \qquad \textrm{on }\partial\Omega,\\
	{A}(u):&= \begin{cases}
	-\Delta u+u^3 \qquad \qquad\qquad\quad\textrm{ }\textrm{ }\textrm{ }\textrm{ in } \Omega_1, \\
	-\Delta_\varepsilon^{p(x)} u \qquad \qquad\qquad\qquad \textrm{ }\textrm{ }\textrm{ in } \Omega_2, \\
	-\textrm{div}\left((\frac{1}{10}+\frac{1}{1+u^2})\nabla u \right) \qquad \textrm{ in } \Omega_3, \\
	-\Delta u \qquad \qquad\qquad\qquad \qquad\textrm{ in } \Omega_4,
	\end{cases}
\end{align*}
where we assume continuity of $u$ across the  interface $I=(-1,1)\times \{1\} \cup \{1\}\times(-1,1)$ (visualized in green in Figure~\ref{fig: config}), $f=10$, $-\Delta_\varepsilon^{p(x)} u:=
-\textrm{div}\left((|\nabla u|^2+\varepsilon^2)^\frac{p(x)-2}{2} \nabla u\right)$, with $p(x)=(1+2x_1)(1+2x_2)$ and $\varepsilon=10^{-10}$. The different partial differential equations on the sub-domains are shown in Figure~\ref{fig: config}. 
The weak form of our model problem is given by: Find $u \in V$ such that
\begin{equation} \label{eq: cont primal problem}
\mathcal{A}(u)(v)=0 \qquad \forall v \in V,
\end{equation}
where 
\begin{align*}
\mathcal{A}(u)(v):=& \int_{\Omega_1} \nabla u \cdot \nabla v  \dx  + \int_{\Omega_1} u^3 v  \dx 
+\int_{\Omega_2} \left((|\nabla u|^2+\varepsilon^2)^\frac{p(x)-2}{2} \nabla u \right) \cdot \nabla v  \dx \\
+&\int_{\Omega_3} \left((\frac{1}{10}+\frac{1}{1+u^2})\nabla u \right) \cdot \nabla v  \dx 
+\int_{\Omega_4} \nabla u \cdot \nabla v  \dx -\int_{\Omega} f v  \dx.
\end{align*}
In this work, we assume that the weak form of  our model problem has a unique solution $u$.
\section{Discretization}
\label{sec: disc}
For the discretization, we decompose $\Omega$ into quadrilateral elements. Here, we use $Q_1$, i.e. bilinear finite elements and $Q_2$, i.e. bi-quadratic finite elements. For the precise definition, we refer the reader to \cite{Endt21,EndtLaRiSchafWi24_book_chapter}. Arising hanging nodes from the adaptive refinement are treated as described in \cite{carey1984finite}. In this work, the discrete spaces with $Q_1$ finite elements and $Q_2$ finite elements,  will be denoted by $V_h$ and $V_h^{(2)}$, respectively. The discretized model problem for $Q_1$ finite elements is given by: Find $u_h \in V_h$ such that 
\begin{equation}\label{eq: discete primal problem}
\mathcal{A}(u_h)(v_h)=0 \qquad \forall v_h \in V_h.
\end{equation}
The discretized model problem  for $Q_2$ finite elements (enriched problem) is given by: Find $u_h^{(2)} \in V^{(2)}$ such that
\begin{equation}\label{eq: enriched primal problem}
\mathcal{A}(u_h^{(2)})(v_h^{(2)})=0 \qquad \forall v_h^{(2)} \in V_h^{(2)}.
\end{equation}
Of course the methodologies in this work can be applied to other finite elements and discretization techniques as well.
\section{Goal Oriented Adaptivity}
\label{sec: Goafem}
In many applications the solution of the problem is not of primary interest, but a certain quantity of interest $J:V\mapsto \mathbb{R}$. In this section, we derive a goal oriented error estimator for the error in the quantity of interest, i.e. $|J(u)-J(\tilde{u})|$ . For goal oriented error estimation, we refer to \cite{KerPruChaLaf2017,EnLaNeWiWo20,BREVIS2021186,BruZhuZwie16,KerPruChaLaf2017,MeiRaVih109,BaRa03,BeRa01,bause2021flexible,bruchhauser2020dual,dolejvsi2021goal,AhEndtSteiWi22,BeuEndtLaWi23,endtmayer2024goal}. This work is based on the dual weighted residual (DWR) method \cite{BeRa01} and an extension to \cite{EndtLaRiSchafWi24_book_chapter,EndtLaWi18_pamm,EndtLaWi20,RanVi2013}.
As in \cite{RanVi2013,EnLaWi18,Endt21}, we introduce an adjoint problem: Find $z \in V$ such that
\begin{equation} \label{eq: cont adjoint problem}
\mathcal{A}'(u)(v,z)=J'(u)(v) \qquad \forall v \in V,
\end{equation}
where $u$ solves \eqref{eq: cont primal problem} and  $\mathcal{A}'$ and $J'$ denotes the Fr\'{e}chet derivatives of $A$ and $J$ with respect to $u$, respectively.
The discretized adjoint problem is given by: Find $z_h \in V_h$ such that
\begin{equation} \label{eq: discrete adjoint problem}
\mathcal{A}'(u_h)(v_h,z_h)=J'(u_h)(v_h) \qquad \forall v_h \in V_h,
\end{equation}
where $u_h$ solves \eqref{eq: discete primal problem}.
Furthermore, the enriched adjoint problem is given by:  Find $z_h^{(2)} \in V_h^{(2)}$ such that
\begin{equation} \label{eq: enriched adjoint problem}
\mathcal{A}'(u_h^{(2)})(v_h^{(2)},z_h^{(2)})=J'(u_h^{(2)})(v_h^{(2)}) \qquad \forall v_h^{(2)} \in V_h^{(2)},
\end{equation}
where $u_h^{(2)}$ solves \eqref{eq: enriched primal problem}.
\begin{theorem}[see \cite{EndtLaRiSchafWi24_book_chapter,EnLaSchaf2023}]\label{thm: Error Representation}
	If $J(u) \in \mathbb{R}$, where $u\in V$ solves our model problem \eqref{eq: cont primal problem}, $u_h^{(2)} \in V_h^{(2)}$ solves the enriched problem \eqref{eq: enriched primal problem} and $z_h^{(2)} \in V_h^{(2)}$ solves the enriched adjoint problem \eqref{eq: enriched adjoint problem}, then for arbitrary but fixed  $\tilde{u}, \tilde{z} \in V_h$
	the error representation formula
	\begin{equation} \label{eq: Errorrepresentation}
	J(u)-J(\tilde{u})= J(u)-J(u_h^{(2)})+ \frac{1}{2}\left(\rho(\tilde{u})(z_h^{(2)}-\tilde{z})+\rho^*(\tilde{u},\tilde{z})(u_h^{(2)}-\tilde{u}) \right)
	-\rho (\tilde{u})(\tilde{z}) + \mathcal{R}^{(3)},
	\end{equation}
	holds, where
	$\rho(\tilde{u})(\cdot) := -\mathcal{A}(\tilde{u})(\cdot)$ ,
	$\rho^*(\tilde{u},\tilde{z})(\cdot) := J'(u)-\mathcal{A}'(\tilde{u})(\cdot,\tilde{z})$ and
	\begin{equation*}
	\begin{split}
	\mathcal{R}^{(3)}:=\frac{1}{2}\int_{0}^{1}&[J'''(\tilde{u}+se_h^{(2)})(e_h^{(2)},e_h^{(2)},e_h^{(2)})
	-\mathcal{A}'''(\tilde{u}+se_h^{(2)})(e_h^{(2)},e_h^{(2)},e_h^{(2)},\tilde{z}+se_h^{*(2)}) \\
	&-3\mathcal{A}''(\tilde{u}+se_h^{(2)})(e_h^{(2)},e_h^{(2)},e_h^{(2)})]s(s-1)\,ds,
	\end{split} 
	\end{equation*}
	with $e_h^{(2)}=u_h^{(2)}-\tilde{u}$ and $e_h^{*(2)} =z_h^{(2)}-\tilde{z}$.
	\begin{proof}
		See \cite{EndtLaRiSchafWi24_book_chapter,EnLaSchaf2023}
	\end{proof}
\end{theorem}
\begin{Remark}
	In this work, we think of $\tilde{u}$ to be an approximation to $u_h$ and $\tilde{z}$ to be an approximation for $z_h$.
\end{Remark}
First of all $J(u)-J(u_h^{(2)})$ is not computable. However, under a saturation assumption \cite{EndtLaRiSchafWi24_book_chapter,EnLaSchaf2023,EndtLaWi21_smart,beuchler2024mathematical}, this part is bounded by the error $|J(u)-J(\tilde{u})|$. 
As discussed in \cite{EndtLaRiSchafWi24_book_chapter} the part $\mathcal{R}^{(3)}$ is of higher order and is often neglected in literature. The part $-\rho (\tilde{u})(\tilde{z})$ is related to the iteration error, which is related to the accuracy of the nonlinear solver. In this work, we solve the nonlinear system until this part is less than $10^{-10}$ and neglect it afterwards. The part $$\eta_h:=\frac{1}{2}\left(\rho(\tilde{u})(z_h^{(2)}-\tilde{z})+\rho^*(\tilde{u},\tilde{z})(u_h^{(2)}-\tilde{u}) \right)$$ is related to the discretization error. In many works the part $\rho(\tilde{u})(z_h^{(2)}-\tilde{z})$ is called the primal part of the discretization error estimator, whereas $\rho^*(\tilde{u},\tilde{z})(u_h^{(2)}-\tilde{u}))$ is called the adjoint part.  As originally done in \cite{RiWi15_dwr}, we use a partition of unity technique to get local error contributions. For more information on the localization we refer to \cite{RiWi15_dwr}.
\section{Multi Goal Oriented Adaptivity}
\label{sec: multigoal}
In many applications, there is more one quantity of interested (QoI). There are several works on multigoal oriented error estimation like \cite{HaHou03,Ha08,PARDO20101953,BruZhuZwie16,EnLaWiWoNeiPAMM2019,KerPruChaLaf2017,EndtWi17,kergrene2018goal,EnLaWi18,EndtLaRiSchafWi24_book_chapter}. Let us assume we are interested in $N \in \mathbb{N}$ quantities of interest (QoIs) $J_1,J_2,\ldots,J_{N}$. In this work, we follow the approach in \cite{EnLaWi18}. We combine the QoIs to a single combined QoI $J_c$ which is given by 
\begin{equation} \label{eq: combined functional}
J_c(v):=\sum_{i=1}^{N}\frac{\textrm{sign}(J_i(u_h^{(2)})-J_i(\tilde{u}))}{|J_i(\tilde{u})|}J_i(v) \qquad \textrm{with } \qquad \textrm{sign}(x):=	
\begin{cases}
\textrm{ }\textrm{ }\textrm{ }1 \qquad x\geq 0, \\
-1 \qquad x < 0.
\end{cases}
\end{equation}
This choice should avoid error cancellation. More information about these error cancellation effects can be found in \cite{Endt21,EndtLaRiSchafWi24_book_chapter,EnLaWi18,EndtWi17}.
Once the quantities of interest are combined to one, we can apply the ideas from Section \ref{sec: Goafem} with $J=J_c$.
\section{Numerical Results}
\label{sec: numerical results}
\begin{figure}[H]
\begin{figure}[H]
	\centering
		\begin{minipage}{0.47\linewidth}
			\includegraphics[width=1.0\linewidth]{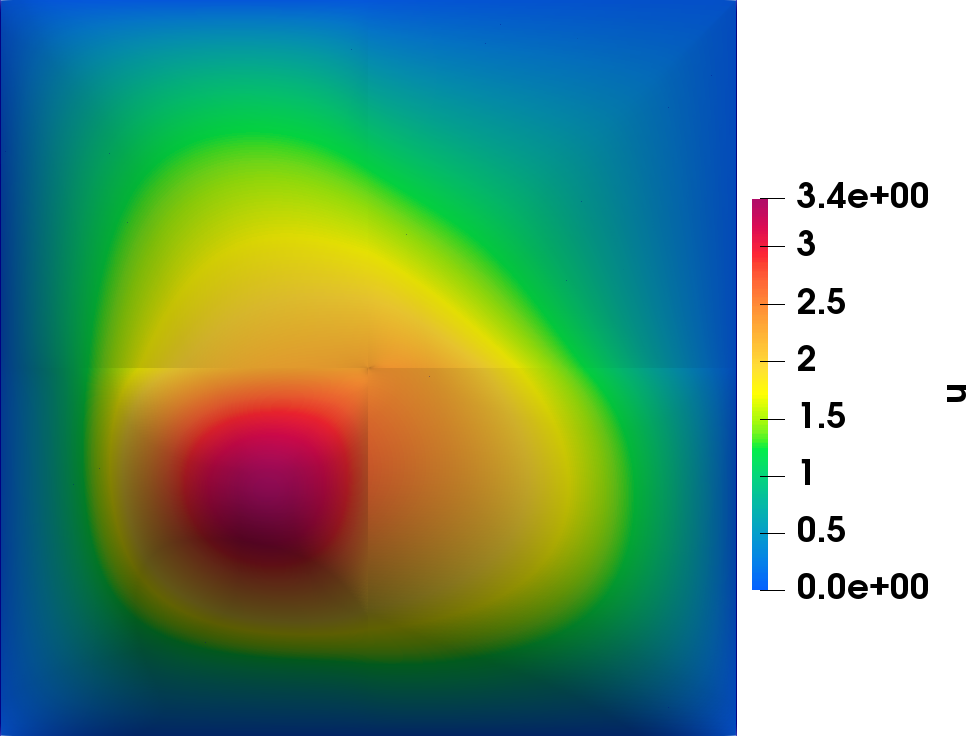}
		\end{minipage}\hfill
		\begin{minipage}{0.47\linewidth}
			\includegraphics[width=0.77\linewidth]{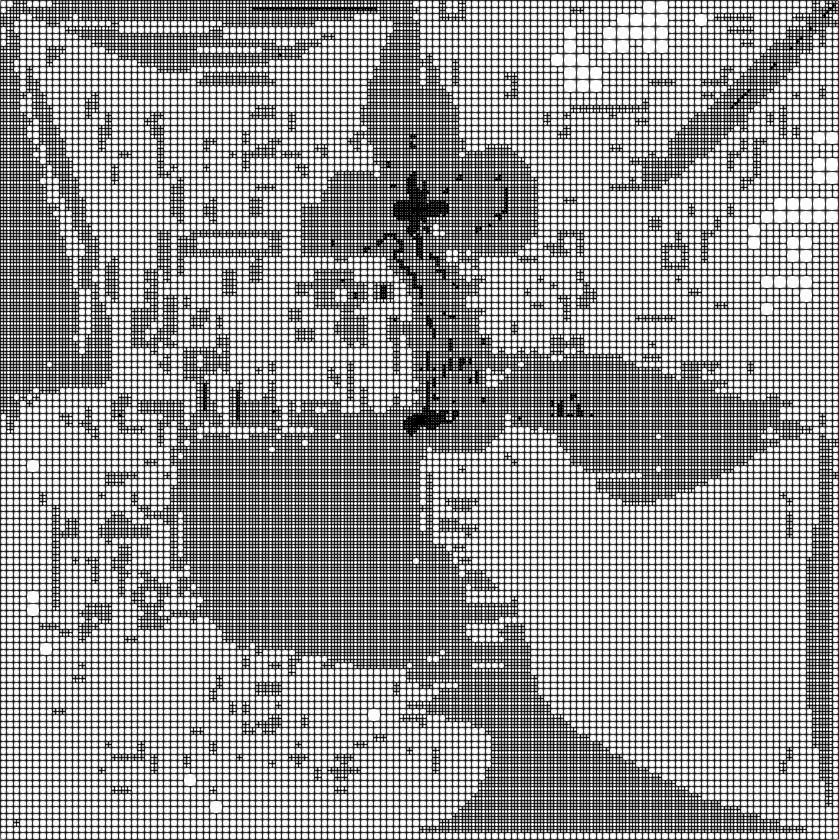}
		\end{minipage}	
		\caption{The solution $u$ (left) and the adaptive mesh after 31 adaptive refinement steps (right).\label{fig: solution + adaptive grid}}
\end{figure}

\begin{figure}[H]
	\centering
	\begin{minipage}{0.47\linewidth}
		\begin{table}[H]
			\begin{tabular}{|c|r|}
				\hline
				Quantity of Interest & \multicolumn{1}{c|}{Reference Value} \\ \hline
				{$J_1(u):=\int_{\Omega_1}u \dx$}            & 0.663674525018                              \\ \hline
				$J_2(u):=\int_{\Omega_2}u \dx$               & 1.011855377549                             \\ \hline
				$J_3(u):=\int_{\Omega_3}u \dx$               & 1.197302269110                    \\ \hline
				$J_4(u):=\int_{\Omega_4}u \dx$               & 1.476574169174                             \\ \hline
				$J_5(u)=u(0,\frac{1}{2})$   &                      1.261144956935                                \\ \hline
			\end{tabular}
		\end{table}
	\end{minipage}\hfill
	\begin{minipage}{0.47\linewidth}
		\includegraphics[width=0.99\linewidth]{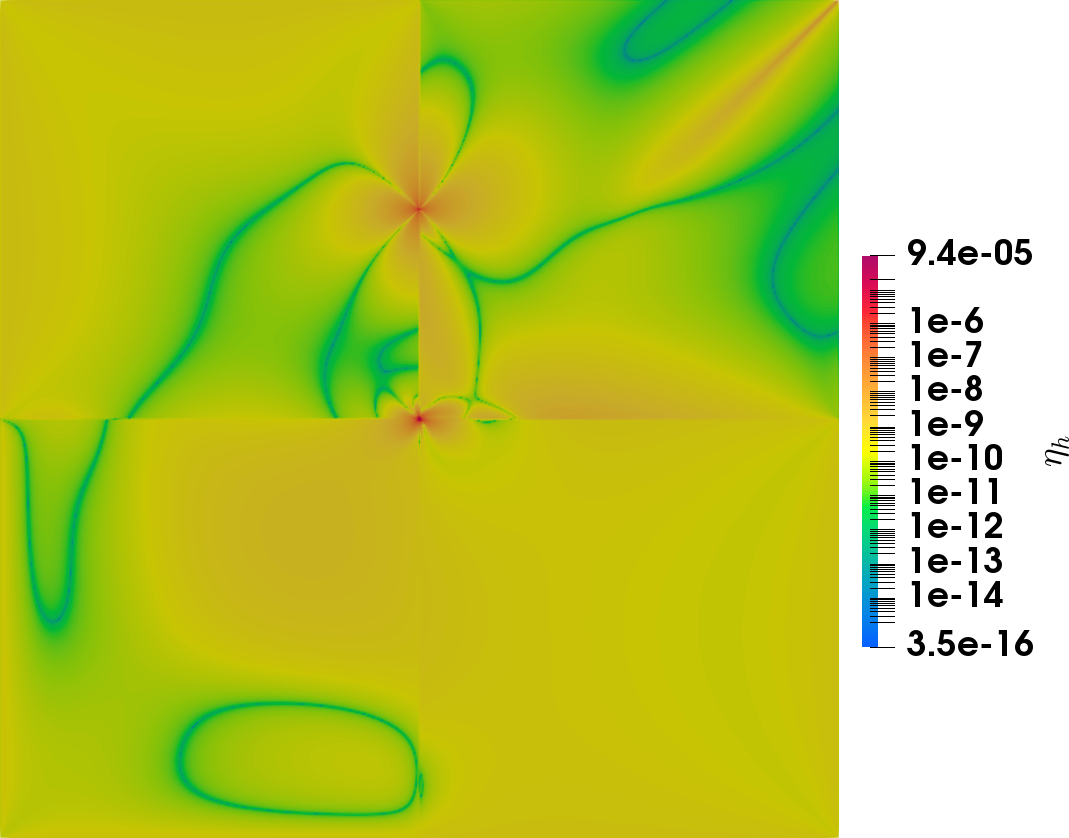}
	\end{minipage}	
	\caption{The quantities of interest in the table (left) and the error distribution resulting from $\eta_h$ with partition of unity localization on an 8 times uniformly refined grid. \label{fig: table + error estimator uniform}.}
\end{figure}

\end{figure}

In the numerical experiments, we used Algorithm~3 from \cite{EndtLaRiSchafWi24_book_chapter} and the finite elements described in Section~\ref{sec: disc}. The initial mesh consists of $4$ elements, namely $\Omega_1,\Omega_2,\Omega_3$ and $\Omega_4$.
The given problem is our model problem described in Section~\ref{sec: model problem}. The solution of our model problem is given in Figure~\ref{fig: solution + adaptive grid} left. The QoIs as well as their reference values are given in Figure~\ref{fig: table + error estimator uniform} left. These QoIs are combined as given in \eqref{eq: combined functional} in Section~\ref{sec: multigoal}.

%

The resulting adaptive refinement after 31 refinement steps is given in Figure~\ref{fig: solution + adaptive grid} (right). There is lot of refinement around the point $x_5=(0,\frac{1}{2})$, where $J_5:=u(x_5)$ is localized. Additionally, there is a lot of refinement around $(0,0)$, where the solution has a singularity. If we compare Figure~\ref{fig: solution + adaptive grid} (right) and Figure~\ref{fig: table + error estimator uniform}, we notice that regions with the largest error contribution on the uniformly refined mesh show the most refinement in the adaptive mesh produced by our algorithm. 
\begin{figure}
	\centering
	\begin{minipage}{0.47\linewidth}
			\includegraphics[width=0.99\linewidth]{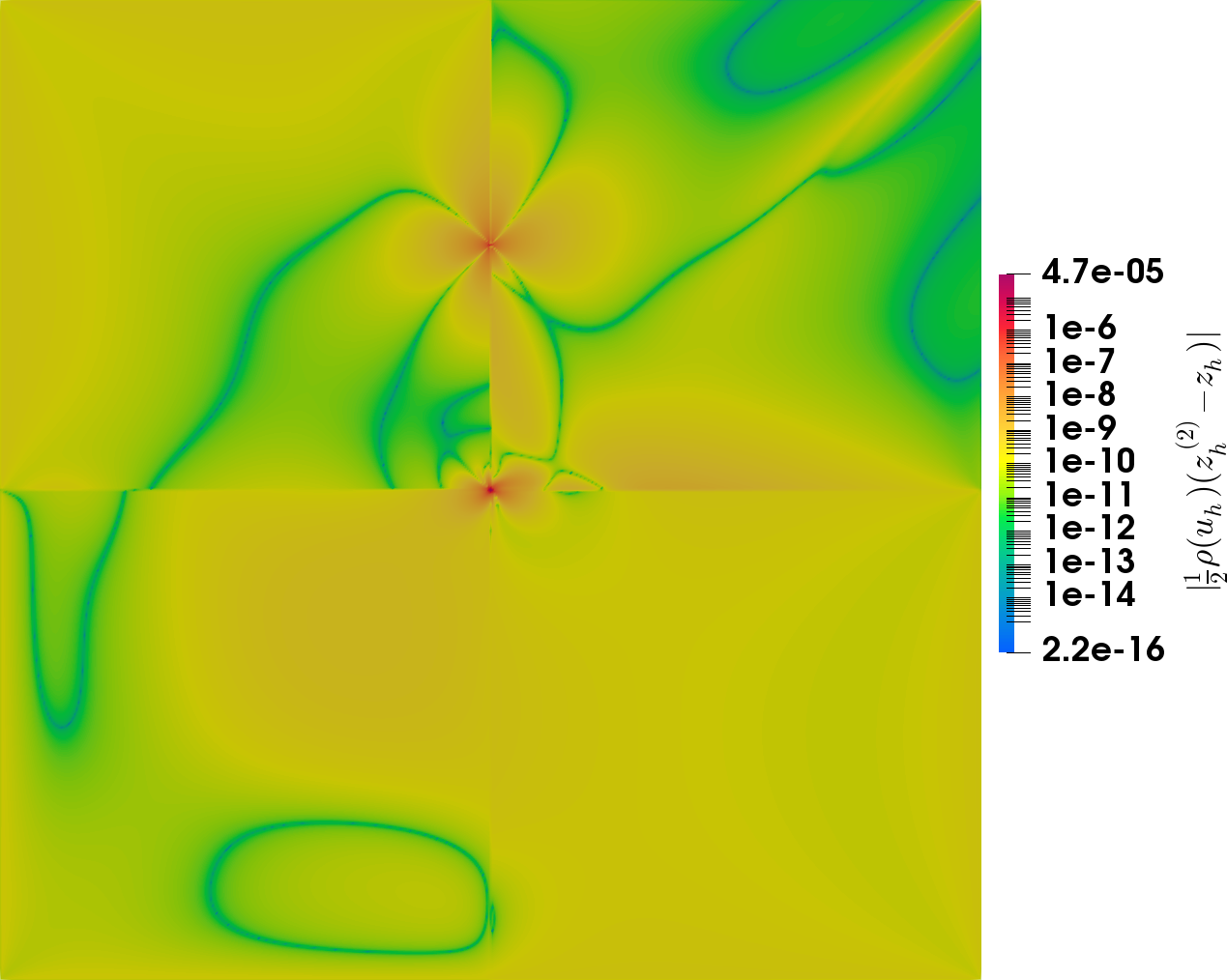}
	\end{minipage}\hfill
	\begin{minipage}{0.47\linewidth}
		\includegraphics[width=0.99\linewidth]{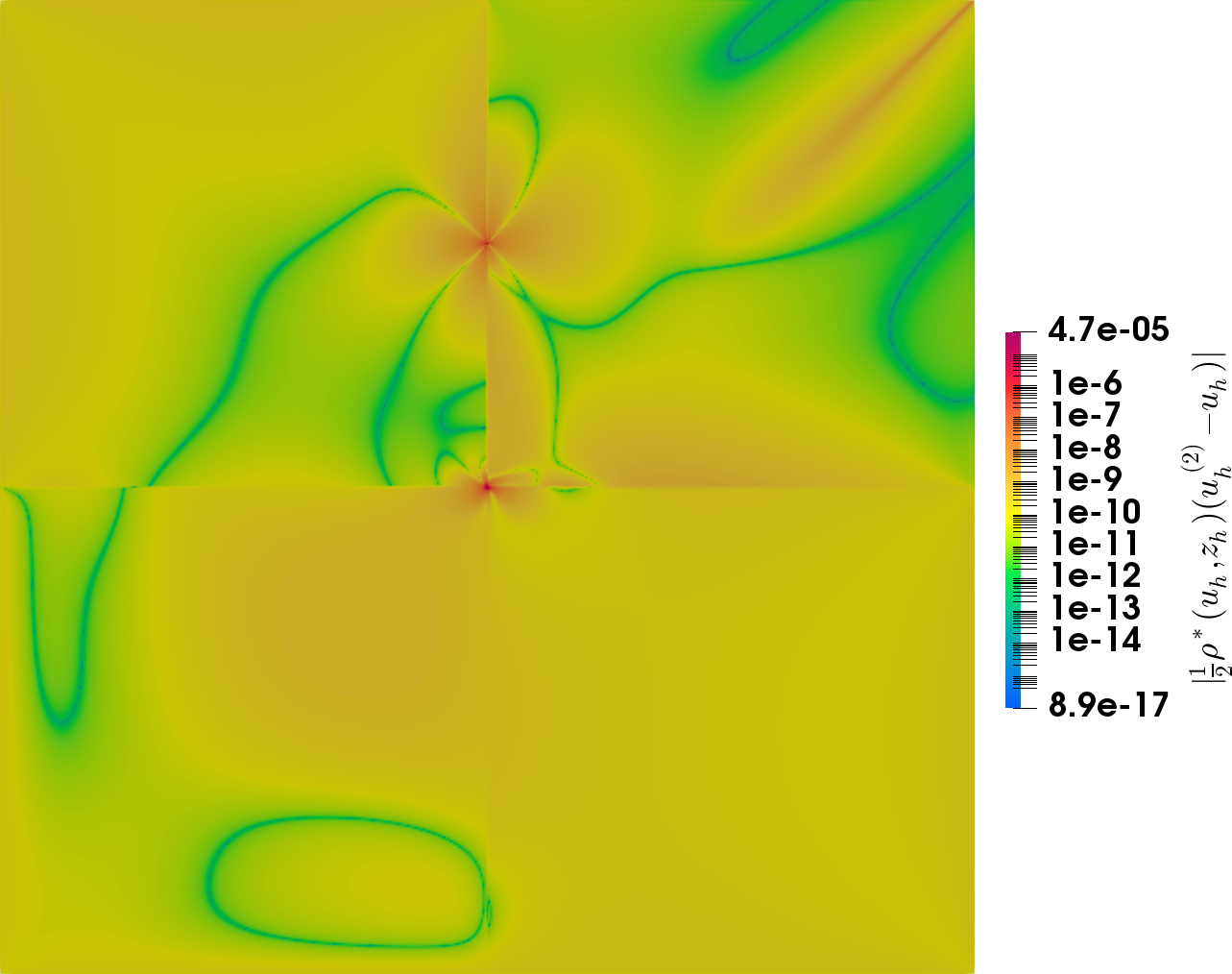}
	\end{minipage}	
	\caption{The localized primal part $\frac{1}{2}\rho(u_h)(z_h^{(2)})$ (left) and adjoint part $\frac{1}{2}\rho^*(u_h,z_h)(u_h^{(2)}-u_h)$ (right) of the error estimator $\eta_h$.\label{fig: localized primal adjoint}}
\end{figure}

In Figure~\ref{fig: localized primal adjoint}, we compare the localized version of the primal and the adjoint part of the error estimator. We observe that both parts almost coincide in the regions $\Omega_1$, $\Omega_3$ and $\Omega_4$. However, in the region $\Omega_2$, where the $p$-Laplace equation is the local PDE, the localizations are more different. This is in accordance with the results in \cite{EndtLaWi20}, where it was shown that for the $p$-Laplace problem the primal and the adjoint part of the the error estimator $\eta_h$ are required.

\begin{figure}


	{	\begin{minipage}{0.47\linewidth}
			\scalebox{0.60}{\input{Figures/Ieff2.tex}} 
			\caption{Effectivity index for $J_{c}$. \label{fig: Ieff}}
		\end{minipage}\hfill
		\begin{minipage}{0.47\linewidth}
			\scalebox{0.60}{\input{Figures/Error2.tex}}
			\caption{Relative errors for $J_1$, $J_2$, $J_3$, $J_4$,  $J_5$ and absolute error for $J_c$. \label{fig: errors}}
		\end{minipage}	
	}	
\end{figure}
In Figure~\ref{fig: Ieff}, the effectivity index, which is the ratio of the error estimator $\eta_h$ and the real error $J(u)-J(u_h)$, is between 0.4 and 1.6. Figure~\ref{fig: errors} shows that the absolute error of $J_c$ bounds all the relative errors of the QoIs $J_1,J_2,J_3,J_4$ and $J_5$. Furthermore, we observe that the error in all quantities of interest decays with rate $O(\mathrm{DoFs}^{-1})$, where $\mathrm{DoFs}$ is the number of degrees of freedom. In Figure~\ref{fig: J3: adaptive+uniform} and Figure~\ref{fig: J5: adaptive+uniform}, the error for uniform refinement has a worse rate than for adaptive refinement. The results for $J_1,J_2$ and $J_4$ are similar to the results shown in Figure~\ref{fig: J3: adaptive+uniform}. 
\begin{figure}[H]

	{	\begin{minipage}{0.47\linewidth}
			\scalebox{0.60}{\input{Figures/ErrorJ22.tex}} 
			\caption{Comparison of the  relative error of $J_3$ for uniform and adaptive refinement. \label{fig: J3: adaptive+uniform}}
		\end{minipage}\hfill
		\begin{minipage}{0.47\linewidth}
			\scalebox{0.60}{\input{Figures/ErrorJ52.tex}}
			\caption{Comparison of the  relative error of $J_5$ for uniform and adaptive refinement. \label{fig: J5: adaptive+uniform}}
		\end{minipage}	
	}	
\end{figure}
\section{Conclusions}
In this work, we applied the dual weighted residual method to a problem with multiple quantities of interest. Furthermore, the model problem itself consist of locally different coercive partial differential equations, which are arranged in a checkerboard pattern. 
The a posteriori goal oriented error estimation for multiple quantities of interest showed a similar convergence rate in all the quantities of interest. Furthermore, the effectivity indices for our model problem are close to 1. Moreover, we observe that the grid generated by the adaptive algorithm shows refinement where the localization of the error estimator on a uniform grid is large. Overall, we observe an excellent performance.

\section*{Acknowledgement}
This work has been supported by the Cluster of Excellence PhoenixD (EXC 2122, Project ID 390833453).
Furthermore, Bernhard Endtmayer greatly acknowledges his current Humboldt Postdoctoral 
Fellowship from the 'Alexander von Humboldt Foundation. Furthermore, i would like to thank Cynthia Schmidt for her social support.

\vspace{\baselineskip}

	\bibliographystyle{abbrv}
\bibliography{lit.bib}

\end{document}